# Approximation Of Logarithm, Factorial And Euler-Mascheroni Constant Using Odd Harmonic Series

**Narinder Kumar Wadhawan**[1*] and **Priyanka Wadhawan**[2]
[1*]Civil Servant, Indian Administrative Service Retired, House No. 563, Panchkula - 134112, Haryana, India. E-mail: narinderkw@gmail.com
[2]Department Of Computer Sciences, Thapar Institute Of Engineering And Technology, Patiala, India, now Program Manager- Space Management (TCS) Walgreen Co. 304 Wilmer Road, Deerfield, Il. 600015 USA, Email: priyanka.wadhawan@gmail.com.

**Abstract**
We have proved in this paper that natural logarithm of consecutive number ratio, x/(x-1) approximates to 2/(2x - 1) where x is a real number except 1. Using this relation, we, then proved, x approximates to double the sum of odd harmonic series having first and last terms 1/3 and 1/(2x - 1) respectively. Thereafter, not limiting to consecutive number ratios, we extended its applicability to all the real numbers. Based on these relations, we, then derived a formula for approximating the value of Factorial x. We could also approximate the value of Euler-Mascheroni constant. In these derivations, we used only and only elementary functions, thus this paper is easily comprehensible to students and scholars.

## 1. Introduction

By applying geometric approach, Leonhard Euler in the year 1748 devised methods of determining natural logarithm of a number [1]. He then extended it to other numbers utilising basic properties of logarithm. Sasaki and Kanada, then worked on determination of precise value of $\log(x)$ using special functions [4]. Different formulae [8] for determining the value of logarithm derived so far, are, 1) $\ln(1 + x) = x - x^2/2 + x^3/3 - \cdots$ up to infinity, where $|x| \leq 1$ and $x \neq 1$, 2) if $Re(x) \geq 1/2$, then $\ln(x) = \sum_{k=2}^{\infty}(x-1)^k/(kx^k)$, 3) $\ln\left(\frac{n+1}{n}\right) = \sum_{k=2}^{\infty}\frac{1}{k(n+1)^k}$, 4) $\ln(x) = \frac{2(x-1)}{(x+1)}\left[\frac{1}{1} + \frac{(x-1)^2}{3(x+1)^2} + \left\{\frac{(x-1)^2}{5(x+1)^2}\right\}^2 + \cdots\right]$. In addition to these series, another alternative to high precision calculation is the formula [8], $ln(x) \simeq \pi/\{2M(1,4/s)\} - m\ln(2)$, where $M$ denotes the arithmetic-geometric mean of 1 and $4/s$, and $s = x2^m > 2^{p/2}$ with $m$ chosen so that $p$ bits of precision is attained. The complexity of computing the natural logarithm (using the arithmetic-geometric mean) is O(*M*(*n*)) ln (*n*). Here *n* is the number of digits of precision at which the natural logarithm is to be evaluated and *M*(*n*) is the computational complexity of multiplying two *n*-digit numbers [8],[10].

Franzen gave the method of approximation of factorial using relation, $\ln(n!) = \sum_{j=1}^{n}\ln(j)$, [2]. Wolfram MathWorld can be referred to for different methods of approximation of factorial derived by different mathematicians [10]. Euler Mascheroni constant appeared for the first time in the paper of Leonhard Euler [7]. Tims and Tyrrell also worked on approximation of this constant [5]. Young gave an inequality for bounding the harmonic number in terms of the hyperbolic cosine for determining this constant [3]. Various methods adopted to approximate this constant find mention in Wikipedia [7] and Wolfram MathWorld [6].

Notwithstanding numerous works already undertaken on calculation of logarithm, factorial of a number and Euler Mascheroni constant, we adopted a completely different, unique and simple approach. Icing on the cake is, it does not involve special functions and that makes it easily comprehensible even to undergraduate students. To start with, a consecutive number ratio $n/(n-1)$ will be expressed in exponential form and then from these ratios, an exponential function will be derived for number $n$. It will be proved that natural logarithm of number $n$ approximates to $2\sum_{x=2}^{n} 1/(2x-1) + 2\sum_{x=2}^{n} 1/\{x^3(2x-1)^2\}$ where symbol $\sum_{x=2}^{n} \ln\{1/(2x-1)\}$ denotes sum of terms $\{1/(2x-1)\}$ where $x$ varies from 2 to $n$. Based on this exponential representation of a number, formula for $n!$ will be derived and value of Euler Mascheroni constant will be approximated. We factorise a number $n$ as shown in Equation (1.1)

$$n = (1+1)(1+1/2)(1+1/3)\ldots\{1+1/(n-1)\} = \prod_{x=2}^{n}\{1+1/(x-1)\} \quad (1.1)$$

where symbol $\prod_{x=2}^{n}\{1+1/(x-1)\}$ denotes product of terms $\{1+1/(x-1)\}$ where $x$ varies from 2 to $n$. Therefore,

$$\ln(n) = \ln(1+1) + \ln\left(1+\frac{1}{2}\right) + \ln\left(1+\frac{1}{3}\right) + \cdots + \ln\left(1+\frac{1}{n-1}\right) = \sum_{x=2}^{n}\ln\left\{1+\frac{1}{x-1}\right\} \quad (1.2)$$

Quantity $\sum_{x=2}^{n} \ln\{1+1/(x-1)\}$ can be approximated to integration of function $f(x)$ with respect to $x$, where $f(x) = \ln\{1+1/(x-1)\}$ and $x$ varies from 2 to $n$. Mathematically,

$$\ln(n) \simeq \int_{2}^{n} \ln\{1+1/(x-1)\}\,dx \quad (1.3)$$

On integration,

$$\ln(n) \simeq \int_{2}^{n} \ln\{1+1/(x-1)\}\,dx \simeq n\cdot\ln(n) - (n-1)\cdot\ln(n-1) - \ln(2).$$

On rearranging,

$$n \simeq (n-1)2^{1/(n-1)}. \quad (1.4)$$

Replacing $x$ with $n$, Equation (1.4) takes the form $x \simeq (x-1)\cdot 2^{1/(x-1)}$. This derivation of representation of $x$ in exponential form proves **Lemma 1.1.**

**Lemma 1.1**: *A number $x$ can be roughly approximated to $(x-1)\cdot 2^{1/(x-1)}$ where $x$ is any positive or negative number.*

However, this approximation suffers serious drawback on account of the fact that at $x=1$, value of $(x-1)\cdot 2^{1/(x-1)}$ is equal to zero, therefore, to obviate this aberration, Equation (1.4) needs correction.

## 2. Theory and Concept

To determine correction, we draw two graphs. First graph is a plot of $\ln\{1+1/(x-1)\}$ taken on Y-axis with variable $x$ taken on X-axis. Area under the plotted curve will correspond to $\int_{2}^{n} \ln\{1+1/(x-1)\}\,dx$. Second graph is a plot of $\ln\{1+1/(x-1)\}$ taken on Y-axis with $x$ taken on X-axis where $x$ varies in steps from 2 to 3, 3 to 4, so on and area under the plotted graph will correspond to $\sum_{x=2}^{n} \ln\{1+1/(x-1)\}$. Kindly refer to 'Figure 1.' Perusal of the graphs reveals that the quantity $\int_{2}^{n} \ln\{1+1/(x-1)\}\,dx$ relates to the area under the smooth curve whereas the quantity $\sum_{x=2}^{n} \ln\{1+1/(x-1)\}$ relates to the area under the step-shaped graph. Since our requirement is the area under the step shaped graph, a correction is necessitated to conform smooth curve to the step shaped graph.

## 2.1 Corrections To Conform Smooth Curve To A Steps Shaped Graph

For conforming area under smooth curve ADGJMP to the area under the step shaped graph ACDFGIJLMOPR, area of triangle ACD is subtracted from magnitude of term 3

$\{T_3 = \ln\left(1 - \frac{1}{3-1}\right)\}$, area of triangle DFG is subtracted from magnitude of term 4 $\{T_4 = \ln\left(1 - \frac{1}{4-1}\right)\}$ so on till the last $xth$ Term. In this way, correction for term 3 $(T_3) = \frac{1}{2}\{\ln\left(1 + \frac{1}{3-1}\right) - \ln\left(1 + \frac{1}{2-1}\right)\}$, correction for term 4 $(T_4) = \frac{1}{2}\{\ln\left(1 + \frac{1}{4-1}\right) - \ln\left(1 + \frac{1}{3-1}\right)\}$, correction for term 5 $(T_5) = \frac{1}{2}\{\ln\left(1 + \frac{1}{5-1}\right) - \ln\left(1 + \frac{1}{4-1}\right)\}$,

... ... ...

and correction for term x $(T_x) = \frac{1}{2}\{\ln\left(1 + \frac{1}{x-1}\right) - \ln\left(1 + \frac{1}{x-1-1}\right)\}$.

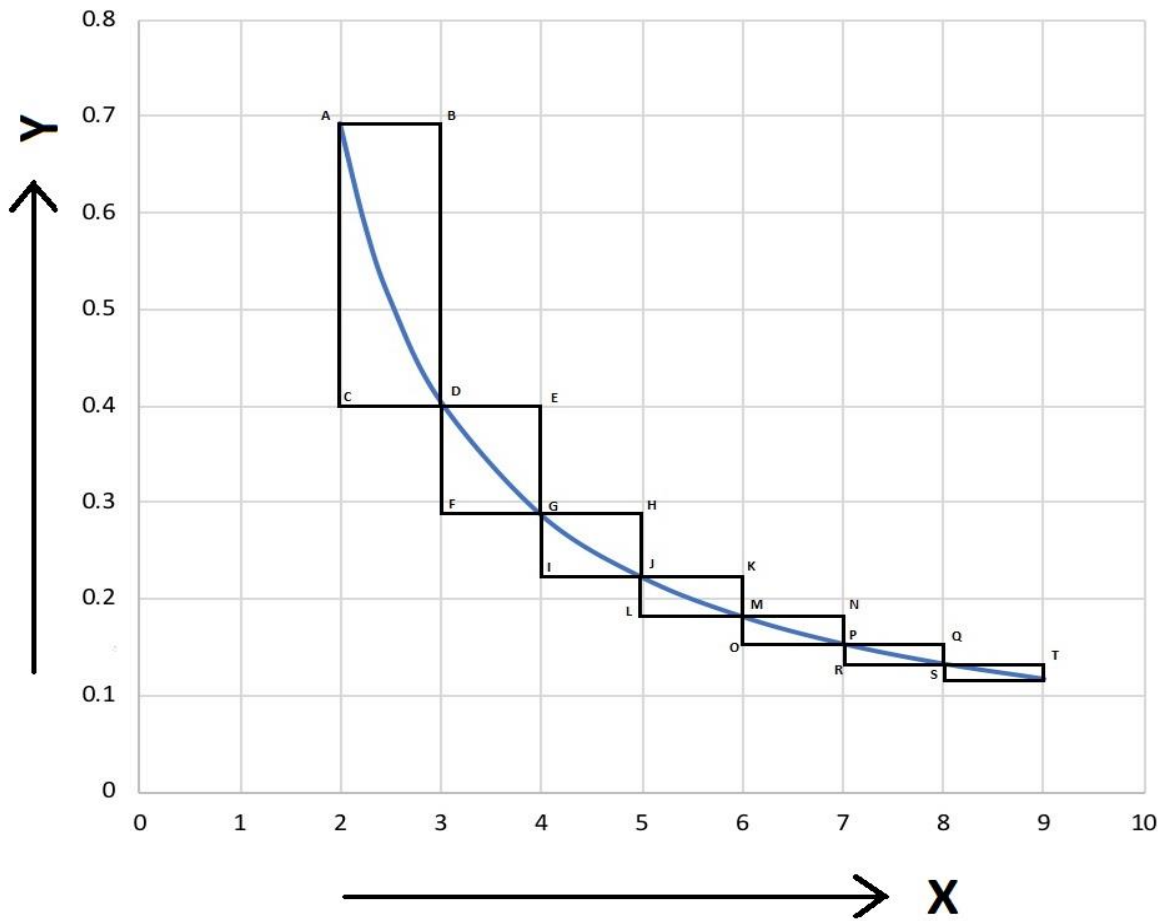

Figure 1 Showing Graphs Of $\ln\left(\frac{x}{x-1}\right)$ With $x$

Being initial condition, magnitude of term 2 $(T_2)$ does not need correction. In this way, the resultant correction is the sum of corrections for all the terms and is equal to $\frac{1}{2}\{\ln\left(1 + \frac{1}{x-1}\right) - \ln(2)\}$. Algebraic addition of this resultant correction to the right hand side of Equation (1.4) yields modified relation,

$$\ln(x) \simeq x \cdot \ln(x) - (x-1) \cdot \ln(x-1) - \ln(2) + \frac{1}{2}\{\ln\left(1 + \frac{1}{x-1}\right) - \ln(2)\}.$$

Or
$$x/(x-1) \simeq 2^{3/(2x-1)}. \tag{2.1}$$

This derivation proves **Lemma 2.1**.

**Lemma 2.1**: *A number x roughly approximates to* $2^{3/(2x-1)}$.

Although Equation (2.1) yields better result, it is not free from error at $x = 1$ or in the vicinity of 1. It is worth mentioning that our assumption in paragraph 2.1 that ACD is a triangle, considering portion A to C a straight line by ignoring the fact, it is a curve, caused error and that error still needs correction.

### 2.1a. Correction Due To Curvature And Also At x → ∞ (Infinity)

In addition to error due to curvature, Equation (2.1) is also not free from error when $x \to \infty$. By definition, $\lim_{x \to \infty} \{x/(x-1)\}^x = e$ where 'e' is Euler's number. Applying Equation (2.1), when $x \to \infty$, quantity $\{x/(x-1)\}^x$ tends to $2^{3/2}$ whereas it should tend to 'e', therefore, there still exists appreciable error warranting additional correction. To eliminate this error when $x \to \infty$, power $3/2$ to the base 2 is replaced by slightly smaller quantity $1/\ln(2)$ and Equation (2.1) gets transformed into

$$x/(x-1) \simeq 2^{2x \cdot \frac{1/\ln(2)}{1-1/2x}}. \qquad (2.2)$$

Since quantity $3/2$ is slightly reduced, therefore, magnitude of denominator $(1 - 1/2x)$ also needs reduction to $\{1 - 1/(2x) - 1/(2x^4)\}$ so as to keep $\{2/\ln(2)\}/\{1 - 1/(2x)\}$ unchained. This quantity $\{1 - 1/(2x) - 1/(2x^4)\}$ was derived by algebraic addition of $-a/x^n$ to $(1 - 1/2x)$ and, then values of 'a' and 'n' were determined by assigning successively different values 1, 2, 3 … to 'a' and 'n' so that quantity $2^{2x \cdot \frac{1/\ln(2)}{1-1/2x-a/x^n}}$ gets equal to $x/(x-1)$. Values of 'a' and 'n' were thus found to be ½ and 4. Since values of $x/(x-1)$ at different values of $x$, are known, therefore, these were utilised to discover values of $a$ and $n$. When values of 'a' and 'n' were determined and since $2^{1/\ln(2)} = e$, therefore, Equation (2.2) gets transformed into

$$x/(x-1) \simeq e^{2/(2x-1-1/x^3)}. \qquad (2.3)$$

Derivation of Equation (2.3) proves **Lemma 2.2.**

**Lemma 2.2**: *A number x approximates to* $(x-1)e^{2/(2x-1-1/x^3)}$ *where x is a positive or negative real number.*

On the basis of equation (2.3), values of $x$ are calculated and are given in **Table 2.1**.

*Table* **2.1**: Approximation Of Numbers Using Equation (2.3)

| $x$ | Value of $x$ Calculated Using Formula (2.3) | Percentage Error | $x$ | Value of $x$ Calculated Using Formula (2.3) | Percentage Error |
|---|---|---|---|---|---|
| 3.5 | 3.493572593 | −.1836402026 | −125 | −125.0000053 | $-4.21428842 \times 10^{-6}$ |
| 5.8 | 5.797266055 | −.0471369806 | −175 | −175.0000027 | $-1.54136032 \times 10^{-6}$ |
| −15.9 | −15.9002932 | −.001844063542 | 750 | 749.9999999 | $-1.97924237 \times 10^{-8}$ |
| −50.1 | −50.1000321 | $-6.41716299 \times 10^{-5}$ | 1500 | 1500 | 0 |
| −100.1 | −100.100008 | $-8.18031141 \times 10^{-6}$ | 2500 | 2500 | 0 |

Perusal of figures given in **Table 2.1** reveals that percentage error decreases with increase in magnitude of the number $x$. For numbers above 1500, percentage error is practically zero.

### 2.1b. Applicability Of Formula (2.3) To Numbers Between −2 To +2

When $x = 3/2$, approximated value of $x$ using Equation (2.3) is 1.540108424 and its percentage error is 2.673894964 which is appreciably large. It can be verified that when $x$ lies between −2 and +2 or is in the vicinity of |2|, percentage error is appreciable and needs

elimination. To overcome this handicap, magnitude of the number $x$ is increased by multiplying it with a large number $m$ and the Equation (2.3) transforms into
$$mx \simeq (mx - 1) \cdot e^{2/\{2mx-1-1/(mx)^3\}}.$$
Division by $m$ restores it to its original value $x$ which is given by relation
$$x \simeq (x - 1/m) \cdot e^{2/\{2mx-1-1/(mx)^3\}}. \tag{2.4}$$
Application of this equation although is compulsorily required for numbers *between* $-2$ *and* $+2$, this can be applied to all numbers for better results. This derivation proves ***Lemma 2.3.***

***Lemma* 2.3:** *A real number $x$ closely approximates to $(x - 1/m) \cdot e^{2/\{2mx-1-1/(mx)^3\}}$ where $m$ is a large integer. Evidently, larger the value of $m$, better will be the result.*

This relation is applicable to all the numbers but is essential for number between $-2$ to $+2$.

**2.1c Examination And Elimination Of Error When $x$ Varies From $+2$ to $-2$**

*Table* **2.2:** Comparison Of Error Using Formula (2.3) And Formula (2.4)

| x | Value of $x$ Calculated Using Formula (4) | Value of x Using Formula (2.4) With $m = 100$ | x | Value of $x$ Calculated Using Formula (2.3) | Value of x Using Formula (2.4) With m as 100 |
|---|---|---|---|---|---|
| 2 | 2.00591 | 1.999999979 | 1.2 | 2.28356 | 1.199999941 |
| 1.8 | 1.82285 | 1.799999974 | 1.0 | 0.00000 | 0.9999999154 |
| 1.6 | 1.66819 | 1.599999967 | 0.7 | $-.13546$ | 0.6999998263 |
| 1.4 | 1.61105 | 1.399999957 | 0.3 | $-.66358$ | 0.2999990268 |

It is explicit from *Table* **2.2** that percentage error is negligible, when the Equation (2.4) was used. It is pertinent to mention that when $x$ is large, magnitude of $1/x^3$ becomes so small that it can be ignored. Equation (2.4), then transforms into
$$x \simeq (x - 1)e^{2/(2x-1)}. \tag{2.5}$$

**2.2. Number Building Blocks And Consecutive Number Ratios**

Quantity $x/(x - 1)$ is a ratio of two consecutive numbers $x$ and $x - 1$ and can be approximated using the Equations (2.3) and (2.4). It is noteworthy that a number $x$ can be built by multiplying consecutive number ratios abbreviated as CNR's. For example, $x = (2/1)(3/2)(4/3) \ldots \{x/(x - 1)\}$. On account of this property, CNR's are also called number building blocks abbreviated as NBB's. Both CNR's and NBB's are synonymous. Data mentioned in the *Table* **2.3** proves the truthfulness of the Equation (2.3).

*Table* **2.3** Comparison Of Values Of Actual CNR's, Logarithm Of CNR's With Calculated Values

| $x$ | Actual CNR $x/(x - 1)$ | Calculated CNR $\frac{2}{e^{2x-1-1/x^3}}$ | Percentage Error | $\ln\{x/(x - 1)\}$ | Calculated $2/(2x - 1 - 1/2x^2)$ | Percentage Error |
|---|---|---|---|---|---|---|
| 1 | $\infty$ | $\infty$ | 0.00000 | $\infty$ | $\infty$ | 0.0000 |
| 2 | 2.00000 | 2.00501 | 0.25081 | 0.69315 | 0.69565 | .36067 |
| 6 | 1.20000 | 1.19949 | 0.04250 | 0.18232 | 0.18189 | .23584 |
| 20 | 1.0526 | 1.05262 | .00110 | .05129 | .05128 | .01949 |
| 100 | 1.0101010 | 1.0101009 | $1.089991 \times 10^{-5}$ | .010050335 | .01005025 | $7.9582 \times 10^{-4}$ |

It is clear from the figures given in the *Table* **2.3** that as $x$ increases, calculated values of CNR's and their logarithms closely approximate to their actual values.

**2.3. Derivation Of Numbers In Exponential Form And Their Logarithms**

By Binomial expansion and ignoring terms having powers more than 3 to the base $x$, $2/(2x - 1 - 1/x^3) = 2/(2x - 1) + 2/\{x^3(2x - 1)^2\}$. Therefore, Equation (2.3) takes the form, $x \simeq (x - 1) \cdot e^{2/(2x-1)} \cdot e^{2/\{x^3(2x-1)^3\}}$ and integers 2, 3, 4, …, $x$ approximate to $[1 \cdot e^{2/3} \cdot e^{2/\{(2^3)(3^2)\}}]$, $[2 \cdot e^{2/5} \cdot e^{2/\{(3^3)(5^2)\}}]$, $[3 \cdot e^{2/7} \cdot e^{2/\{(4^3)(7^2)\}}]$, …, $[(x - 1) \cdot e^{2/(2x-1)} \cdot$

$e^{2/\{(x^3)(2x-1)^2\}}$]. On multiplying integers 2, 3, 4, ..., $x$ and simplifying, we get, $x \simeq e^{2\left(\frac{1}{3}+\frac{1}{5}+\frac{1}{7}+\cdots+\frac{1}{2x-1}\right)} \cdot e^{2\left\{\frac{1}{2^3 \cdot 3^2}+\frac{1}{3^3 \cdot 5^2}+\frac{1}{4^3 \cdot 7^3}+\cdots+\frac{1}{x^3 \cdot (2x-1)^2}\right\}}$. Let the number be denoted by $n$, then

$$n \simeq e^{2\left(\frac{1}{3}+\frac{1}{5}+\frac{1}{7}+\cdots+\frac{1}{2n-1}\right)} \cdot e^{\left\{\frac{1}{2^3 \cdot 3^2}+\frac{1}{3^3 \cdot 5^2}+\frac{1}{4^3 \cdot 7^3}+\cdots+\frac{1}{n^3 \cdot (2n-1)^2}\right\}} \quad (2.6)$$

and

$$\ln(n) \simeq 2\left\{\frac{1}{3}+\frac{1}{5}+\frac{1}{7}+\cdots+\frac{1}{(2n-1)}\right\} + 2\left\{\frac{1}{2^3 \cdot 3^2}+\frac{1}{3^3 \cdot 5^2}+\frac{1}{4^3 \cdot 7^2}+\cdots+\frac{1}{n^3(2n-1)^2}\right\}. \quad (2.7)$$

This proves **Lemma 2.4,**

**Lemma 2.4**: A positive integer $n$ approximates to $e^{2\sum_{x=2}^{n}\left\{\frac{1}{2x-1}+\frac{1}{x^3(2x-1)^3}\right\}}$ and its logarithm $\ln(n)$ to $2\sum_{x=2}^{n}\frac{1}{2x-1} + 2\sum_{x=2}^{n}\frac{1}{x^3(2x-1)^2}$.

Examination of Equation (2.7) reveals that major contributor to approximation of $\ln(n)$ is double the sum of odd harmonic series $2\sum_{x=2}^{n} 1/(2x-1)$. To illustrate the percentage error involved in approximation of $\ln(x)$ using the Equation (2.7), values of actual and calculated $\ln(x)$ are given *Table* **2.4**.

*Table* **2.4**: Comparison Of Actual And Calculated Values of $\ln(x)$

| $x$ | $\ln(x)$ On The Basis Of Equation (2.7) | $\ln(x)$ Actual | Percentage Error | $x$ | $\ln(x)$ On The Basis Of Equation (2.7) | $\ln(x)$ Actual | Percentage Error |
|---|---|---|---|---|---|---|---|
| 2 | 0.69444 | .69315 | .18715 | 10 | 2.29823 | 2.30258 | −.18912 |
| 3 | 1.09740 | 1.09861 | −.10967 | 11 | 2.39347 | 2.39789 | −.18455 |
| 5 | 1.60618 | 1.60944 | −.20242 | 13 | 2.56043 | 2.56495 | −.17611 |
| 7 | 1.94195 | 1.94591 | −.20328 | 15 | 2.70347 | 2.70805 | −.16900 |
| 9 | 2.19296 | 2.19722 | −.19409 | 30 | 3.40119 | 3.49119 | −.13243 |

Perusal of values mentioned in the *Table* **2.4**, reveals that at lower values of $n$, error is comparatively large and at higher value, the error is less.

**2.4a. Derivation Of Formulae For Multiplication And Division Of Two Numbers And Their Corresponding Logarithm**

When $x > y$,

$$x \cdot y \simeq e^{4\left(\frac{1}{3}+\frac{1}{5}+\frac{1}{7}+\cdots+\frac{1}{2y-1}\right)} \cdot e^{2\left(\frac{1}{2y+1}+\frac{1}{2y+3}+\frac{1}{2y+5}+\cdots+\frac{1}{2x-1}\right)} \cdot e^{4\left\{\frac{1}{2^3 \cdot 3^2}+\frac{1}{3^3 \cdot 5^2}+\frac{1}{4^3 \cdot 7^2}+\cdots+\frac{1}{y^3(2y-1)^2}\right\}}$$

$$\cdot e^{2\left\{\frac{1}{(y+1)^3(2y+1)^2}+\frac{1}{(y+2)^3(2y+3)^2}+\frac{1}{(y+3)^3(2y+5)^2}+\cdots+\frac{1}{x^3 \cdot (2x-1)^2}\right\}}, \quad (2.8)$$

$$\ln(x) + \ln(y) \simeq 4\left(\frac{1}{3}+\frac{1}{5}+\frac{1}{7}+\cdots\frac{1}{2y-1}\right) + 2\left(\frac{1}{2y+1}+\frac{1}{2y+3}+\frac{1}{2y+5}+\cdots+\frac{1}{2x-1}\right)$$

$$+4\left\{\frac{1}{2^3 \cdot 3^2}+\frac{1}{3^3 \cdot 5^2}++\frac{1}{4^3 \cdot 7^2}+\cdots+\frac{1}{y^3(2y-1)^3}\right\}$$

$$+2\left\{\frac{1}{(y+1)^3(2y+1)^2}+\frac{1}{(y+2)^3(2y+3)^2}+\frac{1}{(y+3)^3(2y+5)^2}+\cdots+\frac{1}{x^3(2x-1)^2}\right\}, (2.9)$$

$$x/y \simeq e^{2\left(\frac{1}{2y+1}+\frac{1}{2y+3}+\frac{1}{2y+5}+\cdots+\frac{1}{2x-1}\right)} \cdot e^{2\left\{\frac{1}{(y+1)^3 \cdot (2y+1)^2}+\frac{1}{(y+2)^3 \cdot (2y+3)^2}+\frac{1}{(y+3)^3 \cdot (2y+5)^2}+\cdots+\frac{1}{x^3(2x-1)^3}\right\}} \quad (2.10)$$

and

$$\ln(x) - \ln(y) \simeq 2\left(\frac{1}{2y+1}+\frac{1}{2y+3}+\frac{1}{2y+5}+\cdots+\frac{1}{2x-1}\right)$$

$$+2\left\{\frac{1}{(y+1)^3 \cdot (2y+1)^2}+\frac{1}{(y+2)^3 \cdot (2y+3)^2}+\frac{1}{(y+3)^3 \cdot (2y+5)^2}+\cdots+\frac{1}{x^3 \cdot (2x-1)^2}\right\}. \quad (2.11)$$

When $x < y$,

$$x \cdot y \simeq e^{4\left(\frac{1}{3}+\frac{1}{5}+\frac{1}{7}+\cdots+\frac{1}{2x-1}\right)+2\left(\frac{1}{2x+1}+\frac{1}{2x+3}+\frac{1}{2x+5}+\cdots+\frac{1}{2y-1}\right)}$$

$$\cdot e^{4\left\{\frac{1}{2^3\cdot 3^2}+\frac{1}{3^3\cdot 5^2}+\frac{1}{4^3\cdot 7^2}+\cdots+\frac{1}{x^3(2x-1)^2}\right\}+2\left\{\frac{1}{(x+1)^3\cdot(2x+1)^2}+\frac{1}{(x+2)^3\cdot(2x+3)^2}+\cdots+\frac{1}{y^3\cdot(2y-1)^2}\right\}}, \quad (2.12)$$

$$ln(x) + ln(y) \simeq 4\left(\frac{1}{3}+\frac{1}{5}+\frac{1}{7}+\cdots\frac{1}{2x-1}\right) + 2\left(\frac{1}{2x+1}+\frac{1}{2x+3}+\frac{1}{2x+5}\cdots+\frac{1}{2y-1}\right)$$

$$+4\left\{\frac{1}{2^3\cdot 3^2}+\frac{1}{3^3\cdot 5^2}++\frac{1}{4^3\cdot 7^2}+\cdots+\frac{1}{x^3(2x-1)^3}\right\}$$

$$+2\left\{\frac{1}{(x+1)^3(2x+1)^2}+\frac{1}{(x+2)^3(2x+3)^2}+\frac{1}{(x+3)^3(2y+5)^2}+\cdots+\frac{1}{y^3(2y-1)^2}\right\}, \quad (2.13)$$

$$\frac{x}{y} \simeq e^{-2\left(\frac{1}{2x+1}+\frac{1}{2x+3}+\frac{1}{2x+5}+\cdots+\frac{1}{2y-1}\right)} e^{-2\left\{\frac{1}{(x+1)^3\cdot(2x+1)^2}+\frac{1}{(x+2)^3\cdot(2x+3)^2}+\frac{1}{(x+3)^3\cdot(2x+5)^2}+\cdots+\frac{1}{y^3\cdot(2y-1)^2}\right\}}, \quad (2.14)$$

and

$$ln(x) - ln(y) \simeq -2\left(\frac{1}{2x+1}+\frac{1}{2x+3}+\frac{1}{2x+5}\cdots+\frac{1}{2y-1}\right)$$

$$-2\left\{\frac{1}{(x+1)^3\cdot(2x+1)^2}+\frac{1}{(X+2)^3\cdot(2x+3)^2}+\frac{1}{(x+3)^3\cdot(2x+5)^2}+\cdots+\frac{1}{y^3\cdot(2y-1)^2}\right\}. \quad (2.15)$$

### 2.4b. Elimination of Error In Approximation Of Logarithm Of A Number

Let there be a number $p/q$. For elimination of error, we write, $p/q = (m\cdot p)/(m\cdot q)$ where integer $m$ has large value.

**When $p/q > 1$.**

Using equation (2.10), $p/q = (m\cdot p)/(m\cdot q) \simeq$

$$e^{2\left(\frac{1}{2mq+1}+\frac{1}{2mq+3}+\frac{1}{2mq+5}+\cdots+\frac{1}{2mp-1}\right)}e^{2\left\{\frac{1}{(mq+1)^3\cdot(2mq+1)^2}+\frac{1}{(mq+2)^3(2mq+3)^2}+\frac{1}{(mq+3)^3(2mq+5)^2}+\cdots+\frac{1}{(mp)^3(2mp-1)^2}\right\}}, \quad (2.16)$$

and

$$ln(p/q) = \ln\{(m\cdot p)/(m\cdot q)\}$$
$$\simeq 2\left(\frac{1}{2mq+1}+\frac{1}{2mq+3}+\frac{1}{2mq+5}\cdots+\frac{1}{2mp-1}\right)$$
$$+2\left\{\frac{1}{(mq+1)^3(2mq+1)^2}+\frac{1}{(mq+2)^3(2mq+3)^2}+\frac{1}{(mq+3)^3(2mq+5)^2}+\cdots+\frac{1}{(mp)^3(2mp-1)^2}\right\}. \quad (2.17)$$

**When $0 < p/q < 1$,**

$$\frac{p}{q} = \frac{m\cdot p}{m\cdot q} \simeq e^{-2\left(\frac{1}{2mp+1}+\frac{1}{2mp+3}+\frac{1}{2mp+5}+\cdots+\frac{1}{2mq-1}\right)}$$

$$\cdot e^{-2\left\{\frac{1}{(mp+1)^3(2mp+1)^2}+\frac{1}{(mp+2)^3(2mp+3)^2}+\frac{1}{(mp+3)^3(2mp+5)^2}+\cdots+\frac{1}{(mq)^3(2mq-1)^2}\right\}}, \quad (2.18)$$

and

$$ln(p/q) = \ln\{(m\cdot p)/(m\cdot q)\}$$
$$\simeq -2\left(\frac{1}{2mp+1}+\frac{1}{2mp+3}+\frac{1}{2mp+5}+\cdots+\frac{1}{2mq-1}\right)$$
$$-2\left\{\frac{1}{(mp+1)^3(2mp+1)^2}+\frac{1}{(mp+2)^3(2mp+3)^2}+\frac{1}{(mp+3)^3(2mp+5)^2}+\cdots+\frac{1}{(mq)^3(2mq-1)^2}\right\}. \quad (2.19)$$

Since $m$ is appreciably large, values of

$2\left\{\frac{1}{(mq+1)^3(2mq+1)^2}+\frac{1}{(mq+2)^3(2mq+3)^2}+\frac{1}{(mq+3)^3(2mq+5)^2}+\cdots+\frac{1}{(mp)^3(2mp-1)^2}\right\}$ and

$-2\left\{\frac{1}{(mp+1)^3(2mp+1)^2}+\frac{1}{(mp+2)^3(2mp+3)^2}+\frac{1}{(mp+3)^3(2mp+5)^2}+\cdots+\frac{1}{(mq)^3(2mq-1)^2}\right\}$ are

ignorable and Equations (2.16), (2.17), (2.17) and (2.18) transform into Equations (2.20), (2.21), (2.22) and (2.23).

**When $p/q > 1$,** we get

$$p/q = (m\cdot p)/(m\cdot q) \simeq e^{2\left(\frac{1}{2mq+1}+\frac{1}{2mq+3}+\frac{1}{2mq+5}+\cdots+\frac{1}{2mp-1}\right)} \quad (2.20)$$

and

$$\ln\left(\frac{p}{q}\right) = \ln\left(\frac{m\cdot p}{m\cdot q}\right) \simeq 2\left(\frac{1}{2mq+1}+\frac{1}{2mq+3}+\frac{1}{2mq+5}+\cdots+\frac{1}{2mp-1}\right). \quad (2.21)$$

**When $0 < p/q < 1$,**

$$p/q = (m\cdot p)/(m\cdot q) \simeq e^{-2\left(\frac{1}{2mp+1}+\frac{1}{2mp+3}+\frac{1}{2mp+5}+\cdots+\frac{1}{2mq-1}\right)} \quad (2.22)$$

and
$$\ln\left(\frac{p}{q}\right) = \ln\left(\frac{m \cdot p}{m \cdot q}\right) \simeq -2\left(\frac{1}{2mp+1} + \frac{1}{2mp+3} + \frac{1}{2mp+5} + \cdots + \frac{1}{2mq-1}\right). \quad (2.23)$$

Derivations of Equations (2.20), (2.21), (2.22) and (2.23) prove **Lemmas 2.5, 2.6, 2.7 and 2.8**.

*Example*: Let there be $p/q = 1/2$ and assuming $m = 25$, then $(m \cdot p)/(m \cdot q) = 25/50$. On putting these values in equation (2.22), we get, $\ln\left(\frac{1}{2}\right) = \ln\left(\frac{25}{50}\right) \simeq -2\left(\frac{1}{51} + \frac{1}{53} + \frac{1}{55} + \cdots + \frac{1}{99}\right) - 2\left\{\frac{1}{26^3 \cdot 51^2} + \frac{1}{27^3 \cdot 53^2} + \frac{1}{28^3 \cdot 55^2} + \cdots + \frac{1}{50^3 \cdot 99^2}\right\}$. Therefore, $\ln\left(\frac{1}{2}\right) \simeq 2\left(\frac{1}{51} + \frac{1}{53} + \frac{1}{55} + \cdots + \frac{1}{99}\right) = -.693097198$. Actual value of $\ln\left(\frac{1}{2}\right)$ is $-.6931471806$. Percentage error is $0.007210944709$.

Table 2.5 Logarithm Of Numbers Between 0 To 2 After Elimination Of Error

| $\{(m)p\}/\{(m)q\}$ | $\ln(p/q)$ Calculated | $\ln(p/q)$ Actual | Percentage Error |
|---|---|---|---|
| $(25)1/(25)4$ | $-1.38623188$ | $-1.386294361$ | $0.004507060091$ |
| $(25)1/(25)2$ | $-0.693097198$ | $-.6931471806$ | $.00778721778$ |
| $(40)3/(40)4$ | $-0.2876808066$ | $-0.2876820725$ | $4.4001761 \times 10^{-4}$ |
| $(15)9/(15)10$ | $-0.1053600813$ | $-0.1053605157$ | $4.12258637 \times 10^{-4}$ |
| $(30)5/(30)4$ | $0.2231425097$ | $0.2231435513$ | $-4.66791087 \times 10^{-4}$ |
| $(50)3/(50)2$ | $0.4054627934$ | $0.4054651081$ | $-5.70877276 \times 10^{-4}$ |
| $(22)7/(22)4$ | $0.5596121685$ | $0.5596157879$ | $-7.33170063 \times 10^{-4}$ |
| $(10)19/(10)10$ | $0.6418508407$ | $0.6418538862$ | $-4.74480635 \times 10^{-4}$ |

It is explicit from the data given in *Table* 2.5 that percentage error being multiple of $10^{-4}$ is small. By increasing the value of m, the error can be further reduced.

**Lemma 2.5**: A number $p/q$ approximates to $e^{2\left(\frac{1}{2mq+1} + \frac{1}{2mq+3} + \frac{1}{2mq+5} + \cdots + \frac{1}{2mp-1}\right)}$, when $p/q > 1$ and to $e^{-2\left(\frac{1}{2mp+1} + \frac{1}{2mp+3} + \frac{1}{2mp+5} + \cdots + \frac{1}{2mq-1}\right)}$, when $0 < p/q < 1$, where integer m is such that $m \cdot p$ and $m \cdot q$ are both appreciably large. Also, higher the value of $m \cdot p$ and $m \cdot q$, less will be the error.

**Lemma 2.6**: Logarithm of a number $\frac{p}{q}$ approximates to $2\left(\frac{1}{2mq+1} + \frac{1}{2mq+3} + \frac{1}{2mq+5} + \cdots + \frac{1}{2mp-1}\right)$, when $\frac{p}{q} > 1$ and to $-2\left(\frac{1}{2mp+1} + \frac{1}{2mp+3} + \frac{1}{2mp+5} + \cdots + \frac{1}{2mq-1}\right)$, when $0 < \frac{p}{q} < 1$ and integer m is such that $m.p$ and $m.q$ are appreciably large. Also higher the value of $m \cdot p$ and $m \cdot q$, less will be the error.

OR

**Lemma 2.7**: Logarithm of a number $p/q$ approximates to double the sum of odd harmonic series which has terms given by $T = 2 \sum_{x=(q+1/m)}^{p} \{1/(2mx - 1)\}$, where x varies from $(q + 1/m)$ to $p$, $(p/q) > 1$ and integer m is such that $m \cdot p$ and $m \cdot q$ are appreciably large. Also higher the value of $m \cdot p$ and $m \cdot q$, less will be the error.

**Lemma 2.8**: Logarithm of a number $p/q$ approximates to double the sum of odd harmonic series which has terms given by $T = -2 \sum_{x=(p+1/m)}^{q} \{1/(2mx - 1)\}$ where x varies from $(p + 1/m)$ to $q$, $0 < p/q < 1$ and integer m is such that $m \cdot p$ and $m \cdot q$ are appreciably large. Also higher the value of $m \cdot p$ and $m \cdot q$, less will be the error.

**2.5. Algorithm For Determination Of Natural Logarithm Of A Number**
1. Let the given real number be $p/q$.

2. Check if $p/q < 0$. If it is $< 0$, then reject it as there is no logarithm in real quantities for a negative number. Also check if $p/q = 0$ or $p = 0$ or $q = 0$, if either $p$ or $q$ is zero, then go to 3 and reject it as there is no logarithm in real quantities for a number 0 or ∞. When $p/q > 0$, then go to 4.
3. Reject.
4. Put $m = 1, 2, 3, \ldots$ so on till both integers $m \cdot p$ and $m \cdot q$ are larger than say 150. It is observed from the figures mentioned in Table 5, when integers $m \cdot p$ and $m \cdot q$ both are above 100, logarithm of a number has percentage error in multiple of $10^{-4}$. If percentage error less than multiple of $10^{-4}$ is required, value of $m$ needs to be increased so as to make integers $m \cdot p$ and $m \cdot q$, both more than 150. Higher the value of $m \cdot p$ and $m \cdot q$, less will be the error.
5. Record value of $m$ that makes integers $m \cdot p$ and $m \cdot q$, both more than say 150.
6. Check if $p/q > 1$. If it is not, go to 8 otherwise go to 7.
7. Calculate $2\left(\frac{1}{2mq+1} + \frac{1}{2mq+3} + \frac{1}{2mq+5} + \cdots + \frac{1}{2mp-1}\right)$. Let its value be $y$. Go to 9.
8. Calculate $-2\left(\frac{1}{2mp+1} + \frac{1}{2mp+3} + \frac{1}{2mp+5} + \cdots + \frac{1}{2mq-1}\right)$. Let its value be $y$.
9. Print result $y$.
10. Result is $\ln(p/q) \simeq y$.

## 2.6. Factorial Of An Integer 2 And Higher

Factorial of a positive integer is given by relation, $n! = (1)(2)(3)\ldots(n)$ and $\ln(n)$ is given by
$$\ln(n!) = \ln(2) + \ln(3) + \ln(4) + \cdots + \ln(n-2) + \ln(n-1) + \ln(n). \quad (2.24)$$
On substituting values of $\ln 2$, $\ln 3$, $\ln 4, \ldots, \ln n$ obtained using Equation (2.6) in Equation (2.24), we get

$$\ln n! \simeq 2\left(\frac{n-1}{3} + \frac{n-2}{5} + \frac{n-3}{7} + \cdots + \frac{2}{2n-3} + \frac{1}{2n-1}\right)$$
$$+ 2\left\{\frac{n-1}{2^3 \cdot 3^2} + \frac{n-2}{3^3 \cdot 5^2} + \frac{n-3}{4^3 \cdot 7^2} + \cdots + \frac{n-(n-2)}{(n-1)^3(2n-3)^2} + \frac{n-(n-1)}{n^3 \cdot (2n-1)^2}\right\}.$$

For better clarity, step by step simplification is given below. On rearranging terms,

$$\ln n! \simeq 2n\left(\frac{1}{3} + \frac{1}{5} + \frac{1}{7} + \cdots + \frac{1}{2n-3} + \frac{1}{2n-1}\right)$$
$$- \left(1 - \frac{1}{3} + 1 - \frac{1}{5} + 1 - \frac{1}{7} + \cdots + 1 - \frac{1}{2n-3} + 1 - \frac{1}{2n-1}\right)$$
$$+ 2\left\{\frac{n}{2^3 \cdot 3^2} + \frac{n}{3^3 \cdot 5^2} + \frac{n}{4^3 \cdot 7^2} + \cdots + \frac{n}{(n-1)^3(2n-3)^2} + \frac{n}{n^3(2n-1)^2}\right\}$$
$$- 2\left\{\frac{1}{2^3 \cdot 3^2} + \frac{2}{3^3 \cdot 5^2} + \frac{3}{4^3 \cdot 7^2} + \cdots + \frac{(n-2)}{(n-1)^3(2n-3)^2} + \frac{(n-1)}{n^3(2n-1)^2}\right\}.$$

Or

$$\ln n! \simeq n \cdot \ln(n) + \left(\frac{1}{3} + \frac{1}{5} + \frac{1}{7} + \cdots + \frac{1}{2n-3} + \frac{1}{2n-1}\right) - (n-1)$$
$$- \left[\left(\frac{1}{2^3 \cdot 3} - \frac{1}{2^3 \cdot 3^2}\right) + \left(\frac{1}{3^3 \cdot 5} - \frac{1}{3^3 \cdot 5^2}\right) + \left(\frac{1}{4^3 \cdot 7} - \frac{1}{4^3 \cdot 7^2}\right) + \cdots \right.$$
$$\left. + \left\{\frac{1}{n^3 \cdot (2n-1)} - \frac{1}{n^3 \cdot (2n-1)^2}\right\}\right].$$

Or

$$\ln n! \simeq \left(n + \frac{1}{2}\right)\ln(n) - (n-1)$$
$$- \left\{\left(\frac{1}{2^3 \cdot 3}\right) + \left(\frac{1}{3^3 \cdot 5}\right) + \left(\frac{1}{4^3 \cdot 7}\right) + \cdots + \frac{1}{(n-1)^3(2n-3)} + \frac{1}{n^3(2n-1)}\right\}.$$

Or

$$\ln n! \simeq \left(n + \frac{1}{2}\right)\ln(n) - (n-1) - \sum_{x=2}^{n}\left\{\frac{1}{x^3(2x-1)}\right\}. \quad (2.25)$$

On decomposing $\frac{1}{x^3(2x-1)}$ into partial fractions, we get, $\sum [1/\{x^3(2x-1)\}] = \sum\{-1/x^3 - 2/x^2 - 4/x + 8/(2x-1)\}$. Since $\sum\{-1/x^3 - 2/x^2 - 4/x + 8/(2x-1)\}$ approximates to $\int\{-1/x^3 - 2/x^2 - 4/x + 8/(2x-1)\}dx$, therefore,

$$\sum\{-1/x^3 - 2/x^2 - 4/x + 8/(2x-1)\} \simeq \int\{-1/x^3 - 2/x^2 - 4/x + 8/(2x-1)\}dx$$

Or

$$\sum[1/\{x^3(2x-1)\}] \simeq 1/(2x^2) + 2/x - 4 \cdot \ln(x) + 4 \cdot \ln(2x-1) + C.$$

At $x = 2$, $\sum[1/\{x^3(2x-1)\}]$ is equal to $1/\{3(2^3)\}$. Therefore, at $x = n$, we obtain,

$$\sum\left\{\frac{1}{x^3(2x-1)}\right\} \simeq \frac{1}{2^3 \cdot 3} + \frac{1}{2}\left(\frac{1}{n^2} - \frac{1}{2^2}\right) + 2\left(\frac{1}{n} - \frac{1}{2}\right) + 4 \cdot \ln\left(1 - \frac{1}{2n}\right) - 4 \cdot \ln(3) + 8 \cdot \ln(2).$$

Or

$$\sum_{x=2}^{n}\left\{\frac{1}{x^3(2x-1)}\right\} \simeq 0.06739495647 + 2\left(\frac{1}{n} + \frac{1}{4n^2}\right) + 4 \cdot \ln\left(1 - \frac{1}{2n}\right).$$

On putting this value $\sum_{x=2}^{n}[1/\{x^3(2x-1)\}]$ in equation (2.25),

$$\ln n! \simeq \left(n + \frac{1}{2}\right) \cdot \ln(n) - (n-1) - 0.06739495647 - 2\left(\frac{1}{n} + \frac{1}{4n^2}\right) - 4 \cdot \ln\left(1 - \frac{1}{2n}\right).$$

Or

$$n! \simeq e^{.9326050435} \cdot n^{\frac{1}{2}} \cdot (n/e)^n \cdot e^{-2\left(\frac{1}{n} + \frac{1}{4n^2}\right)} \cdot \{1 - 1/(2n)\}^{-4}. \quad (2.26)$$

On account of our assumption made in paragraph 2.3 that term $1/\{x^6(2x-1)^3\}$ and other terms containing higher power of $x$ are ignored and also another assumption that $\sum\{-1/x^3 - 2/x^2 - 4/x + 8/(2x-1)\}$ approximates to $\int\{-1/x^3 - 2/x^2 - 4/x + 8/(2x-1)\}dx$, error has crept in the approximation of $n!$. However, this error is appreciably reduced if equation (2.26) is modified to

$$n! \simeq \sqrt{e^{1.83788} \cdot n} \cdot (n/e)^n \cdot e^{-2\left(\frac{1}{n} + \frac{10}{33.n^2}\right)} \cdot (1 - 200/387n)^{-4}. \quad (2.27)$$

Based on equation (2.27), factorial of some integers are calculated and are given in the **Table 2.6**. Perusal of data mentioned in the Table, reveals that maximum percentage error is approximately 0.5 and this gets reduced as $n$ increases. For 160!, it is approximately 0.01. Further, this derivation proves the purpose that '*summation of odd harmonic series facilitates approximation of factorial of a number.*'

*Table 2.6:* Percentage Error Associated With Approximation Of Factorial Using Equation (2.27)

| $n!$ | Actual $n!$ | Calculated $n!$ | Percentage Error | $n!$ | Actual $n!$ | Calculated $n!$ | Percentage Error |
|---|---|---|---|---|---|---|---|
| 2! | 2 | 2.00584 | .29198 | 45! | $1.19622221 \times 10^{56}$ | $1.19575474 \times 10^{56}$ | −.039078 |
| 3! | 6 | 5.96749 | −.54182 | 60! | $8.31860099 \times 10^{81}$ | $8.32098711 \times 10^{81}$ | −.02867 |
| 4! | 24 | 23.87311 | −.52869 | 75! | $2.48091408 \times 10^{109}$ | $2.48035295 \times 10^{109}$ | −.02262 |
| 5! | 120 | 119.46289 | −.44759 | 95! | $1.03299785 \times 10^{148}$ | $1.03281577 \times 10^{148}$ | −.01763 |

| 10 | 362880 | 3621048 | −.21360 | 110! | $1.58824554 \times 10^{178}$ | $1.58800549 \times 10^{178}$ | −.01511 |
| --- | --- | --- | --- | --- | --- | --- | --- |
| 15! | $1.30767 \times 10^{12}$ | $1.305926 \times 10^{12}$ | −.13371 | 125! | $1.88267718 \times 10^{209}$ | $1.88242823 \times 10^{209}$ | −.01322 |
| 25! | $1.54996 \times 10^{25}$ | $1.55112 \times 10^{25}$ | -.074715 | 140! | $1.34620125 \times 10^{241}$ | $1.34604309 \times 10^{241}$ | −.01175 |
| 35! | $1.03331 \times 10^{40}$ | $1.03278 \times 10^{40}$ | −.05141 | 160! | $4.71472364 \times 10^{284}$ | $4.71424166 \times 10^{284}$ | −.01022 |

## 2.7. Number Constant

Quantity $2 \sum_{x=2}^{n}[1/\{x^3(2x-1)^2$ appearing in the formula for approximation of logarithm of a number $n$ is equal to *Number Constant* provided $n \to \infty$ and for calculating its numerical value, quantity $\frac{1}{x^3(2x-1)^2}$ is decomposed into partial fractions $\frac{12}{x} - \frac{24}{2x-1} + \frac{4}{x^2} + \frac{1}{x^3} + \frac{8}{(2x-1)^2}$, therefore, $\sum \frac{1}{x^3(2x-1)^2} \simeq \int \left\{\frac{12}{x} - \frac{24}{2x-1} + \frac{4}{x^2} + \frac{1}{x^3} + \frac{8}{(2x-1)^2}\right\} dx$ where $x$ varies from 2 to $n$.

On integration and simplification,

$$2 \sum_{x=2}^{n} \frac{1}{x^3(2x-1)^2} \simeq 24 \cdot \ln\left(\frac{x}{2x-1}\right) - \frac{8}{x} - \frac{1}{x^2} - \frac{8}{2x-1} + 16.67560703904,$$

and

$$\lim_{n\to\infty} 2 \sum_{x=2}^{n} \frac{1}{x^3(2x-1)^2} \simeq -24 \cdot \ln(2) + 16.67560703904 \simeq .040074705601703. \quad (2.28)$$

According to Equation (2.7), $\ln(n) - 2 \sum_{x=2}^{n} \frac{1}{(2x-1)} \simeq 2 \sum_{x=2}^{n} \frac{1}{x^3(2x-1)^2}$, therefore, when $n \to \infty$,

$$\ln(n) - 2 \sum_{x=2}^{n} \frac{1}{(2x-1)} \simeq 2 \sum_{x=2}^{n} \frac{1}{x^3(2x-1)^2} \simeq N_r \simeq 0.040074705601703. \quad (2.29)$$

where $N_r$ is a constant called *Number Constant*.

## 2.8. Approximation Of Euler-Mascheroni Constant $\gamma$

According to equation (2.29), when $n \to \infty$, $\ln(n) \simeq 2 \sum_{x=2}^{n} 1/(2x-1) + N_r$. On adding and subtracting $2 - \ln(2)$, the Equation (2.29) can be rewritten as

$$\lim_{n\to\infty} \ln(n) = \lim_{n\to\infty} 2 \left(\frac{1}{3} + \frac{1}{5} + \frac{1}{7} + \cdots \frac{1}{2n-1}\right) + \{2 - \ln(2)\} - \{2 - \ln(2)\} + N_r.$$

Or

$$\lim_{n\to\infty} \ln(n) = \lim_{n\to\infty} 2 \left(\frac{1}{3} + \frac{1}{5} + \frac{1}{7} + \cdots \frac{1}{2n-1}\right) + \left\{2 - \left(1 - \frac{1}{2} + \frac{1}{3} - \frac{1}{4} + \cdots upto \infty\right)\right\} - \{2 - \ln(2)\} + N_r,$$

since $\ln(2) = \left(1 - \frac{1}{2} + \frac{1}{3} - \frac{1}{4} + \cdots up\ to\ \infty\right)$. On adding $\ln(2)$ in RHS and LHS, rearranging and simplifying,

$$\lim_{n\to\infty} \left(1 + \frac{1}{2} + \frac{1}{3} + \frac{1}{5} + \frac{1}{7} + \cdots \frac{1}{2n-1}\right) - \lim_{n\to\infty} \ln(n) - \ln(2) \simeq 2 - 2\cdot\ln(2) - N_r.$$

Or

$$\lim_{n\to\infty} \sum_{x=1}^{2n} \frac{1}{x} - \lim_{n\to\infty} \ln(2n) \simeq 2 - 2 \cdot \ln(2) - N_r.$$

Let $p = 2n$. When $n \to \infty$, $p$ also tends to $\infty$ and

$$\lim_{p\to\infty} \sum_{x=1}^{p} (1/x) - \lim_{p\to\infty} \ln(p) \simeq 2 - 2 \cdot \ln(2) - N_r. \quad (2.30)$$

By definition, Euler-Mascheroni Constant $\gamma = \lim_{p\to\infty} \sum_{x=1}^{p}(1/x) - \lim_{p\to\infty} \ln(p)$. $\quad (2.31)$

Therefore,     Euler $-$ Mascheroni Constant $\simeq 2 - 2 \cdot \ln(2) - N_r$. $\quad (2.32)$

On substituting the value of $N_r$,

$$\text{Euler} - \text{Mascheroni Constant} \simeq 2 - 2 \cdot \ln(2) - \sum_{2}^{\infty} \left\{ \frac{2}{n^3 \cdot (2n-1)^2} \right\}. \qquad (2.33)$$

**2.8a. Approximation Of Numerical Value Of Euler Mascheroni Constant $\gamma$**

On putting the value of $N_r$ given by equation (2.29) in equation (2.30), Euler Mascheroni Constant $\simeq 0.5736309333$. Actual value of Euler Mascheroni Constant [7] $\gamma$ is $0.577215664901$. There is an error of about $-0.6$ percent which is attributed to the fact that exponential quantity $2/(2n-1-1/n^3)$ appearing in equation (2.3), is approximated to $2/(2n-1) + 2/\{n^3(2n-1)^2\}$ ignoring term $2/\{n^6(2n-1)^3\}$ and other terms containing higher powers of $n$. Also in paragraph 2.7, $\lim_{n \to \infty} \sum_{x=2}^{n} 1/\{x^3(2x-1)^2\}$ is approximated to $\int_{2}^{\infty} \{12/x - 24/(2x-1) + 4/x^2 + 1/x^3 + 8/(2x-1)^2\} dx$. Even Equation (2.1) derived for exponential representation of numbers gives approximate results. All these factors accumulated the error to $-0.6$ p.c.

**3. Conclusions And Results**

Overview of what has been described in the paper, makes it amply clear that the golden key to open the doors to approximation of logarithm, factorial of a number or Euler Mascheroni constant, is the representation of a number in exponential form and that was made possible by factorisation of the number. In general, when a number $x$ is factorised, its $nth$ factor is given by the relation $T_n = n/(n-1)$ which is simply a ratio of two consecutive integers $n$ and $(n-1)$ but is extraordinary enough to build a number. It will not be out of context to state that it bears a strong analogy with human cells. Human cells trillions in number compose the body and NBB's, in the same way, though limited to $(x-1)$ in number, compose a number $x$. To illustrate how NBB's compose a number, we present the relation between $x$ and NBB's.

$$x = (2/1) \cdot (3/2) \cdot (4/3) \ldots \{x/(x-1)\}.$$

It is obvious, NBB's $(2/1), (3/2), (4/3), \ldots, \{x/(x-1)\}$ are $(x-1)$ in number and when multiplied generate a number $x$. In normal practice, a number $x$ is envisaged as that which has magnitude equal to what we get when 1 is added $x$ times. This concept of numbers by addition, arises on account of the fact, we are taught mathematics starting with 'counting of the numbers,' in kindergarten and basing thereupon, we distinguish one number from the other on account of weight acquired on accumulation of *unities* in it. We do addition and subtraction corresponding integers to our fingers that is why fingers are called digits. With these strong prejudices, we are unable to envisage an integer as product of numbers. Thinking out of box, we have considered, in this paper, an integer to be product of NBB's. Based upon that we give some examples. $(101/3)$ is product of BBB's $(4/3), (5/4), (7/8) \ldots, (101/100)$ and 20 is product of NBB's $(4/3), (4/3), (5/4), (5/4), (7/6), \ldots, (36/35)$. Using NBB's, a number say $x$ can be represented in infinite ways. This is the crux of the research highlighted in the paper.

Once NBB's are known, a function to approximate NBB's is devised using the identity $\ln(n) = \sum_{2}^{n} x/(x-1)$. Approximating $\sum_{2}^{n} x/(x-1)$ to $\int_{2}^{x} \ln\{x/(x-1)\} dx$, we derived a relation, $x \simeq e^{2/(2x-1-1/x^3)}$. This relation is, then simplified to $x \simeq e^{2/(2x-1)} \cdot e^{2/\{x^3(2x-1)^2\}}$ for its easy applicability to approximation of a number, a factorial and Euler Mascheroni constant. On the basis of this relation, we derived, $\ln(x) \simeq 2 \sum_{n=2}^{x} 1/(2n-1) + 2 \sum_{n=2}^{x} [1/\{n^3(2n-1)^2\}]$. Right hand side of this equation has terms $2 \sum_{n=2}^{x} [1/\{x^3(2x-1)^2\}]$, which are difficult to calculate, therefore, a method is innovated to get rid of these. It is observed, $\frac{1}{2^3 \cdot 3^2} > \frac{1}{3^3 \cdot 5^2} > \frac{1}{4^3 \cdot 7^2} > \cdots > \frac{1}{x^3(2x-1)^2}$, therefore, initial terms like $\frac{1}{2^3 \cdot 3^2}, \frac{1}{3^3 \cdot 5^2}, \frac{1}{4^3 \cdot 7^2}, \ldots$ have substantive values. To reduce the effect of these terms, multiplier $m$ having large value, is used so as to make $x$ as $(m \cdot x)/m$. In this way, the equations are modified to $x \simeq e^{2 \sum_{n=m+1}^{mx} \{1/(2n-1)\}}$.

$e^2 \sum_{n=m+1}^{mx} [1/\{n^3(2n-1)^2\}]$ and $ln(x) \simeq 2 \sum_{n=m+1}^{mx} 1/(2n-1) + 2 \sum_{n=m+1}^{mx}[1/\{n^3(2n-1)^2\}]$. Since value of $m$ is appreciably large, therefore, sum of series $2\left\{\frac{1}{(m+1)^3 \cdot (2m+1)^2} + \frac{1}{(m+2)^3 \cdot (2m+3)^2} + \frac{1}{(m+3)^3 \cdot (2m+5)^2} + \cdots + \frac{1}{(mx)^3 \cdot (2mx-1)^2}\right\}$ can be ignored. In that situation, $x \simeq e^2 \sum_{n=m+1}^{mx} \{1/(2n-1)\}$ and $ln(x) \simeq 2 \sum_{n=m+1}^{mx} 1/(2n-1)$. *It is explicit from above equations that double the sum of odd harmonic series* $2\sum_{n=m+1}^{mx} 1/(2n-1)$ *approximates to natural logarithm of* $x$. These equations were then, applied to derive formulae for product and division of two numbers.

*Factorial of a positive integer* :

By definition, $ln(n!) = ln(2) + ln(3) + ln(4), + \cdots + ln(n)$. Substituting values of $ln(2), ln(3), ln(4), \ldots, ln(n)$ as given by Equation (2.7) in above equation and, then simplifying, we obtain, $n! \simeq e^{.9326050435} \cdot n^{\frac{1}{2}} \cdot (n/e)^n \cdot e^{-2\left(\frac{1}{n} + \frac{1}{4n^2}\right)} \cdot \{1 - 1/(2n)\}^{-4}$ and after correction,

$$n! \simeq \sqrt{e^{1.83788} \cdot n} \cdot (n/e)^n \cdot e^{-2\left(\frac{1}{n} + \frac{10}{33.n^2}\right)} \cdot (1 - 200/387n)^{-4}$$

*Euler Mascheroni. Constant:*

By definition, Euler Mascheroni. constant is equal to $\lim_{p \to \infty} \sum_{x=1}^{p}(1/x) - \lim_{p \to \infty} ln(p)$. Using Equation (2.7) and simplifying, we obtain equation, $Euler - Mascheroni\ Constant \simeq 2 - 2 \cdot ln(2) - \sum_{n=2}^{\infty}[2/\{n^3(2n-1)^2\}]$. Its value, on calculation, is found to approximate to 0.577215664901.

On summing up, derivation of equation, $ln(x) \simeq 2 \sum_{n=2}^{x} 1/(2n-1)$, facilitated approximation of $n!$ and Euler Mascheroni constant. In other words, double the sum of odd harmonic series approximates $ln(x)$ and facilitates approximation of $x!$ and γ. This is what the title of the research says.

### 4. Acknowledgements

We acknowledge the help provided by the wonderful website https://www.desmos.com in calculating the values of tedious exponential terms and successive multiplications of terms.